\documentclass[11pt]{amsart}

\usepackage{amssymb}
\usepackage{natbib}
\usepackage{amscd}

\theoremstyle{definition}
\newtheorem{definition}{Definition}[section]
\newtheorem{example}[definition]{Example}

\theoremstyle{plain}
\newtheorem{theorem}[definition]{Theorem}
\newtheorem{lemma}[definition]{Lemma}
\newtheorem{corollary}[definition]{Corollary}
\newtheorem{proposition}[definition]{Proposition}

\newcommand{\df}[2]{\delta_{#2}#1}
\newcommand{\cdeg}[1]{\mathrm{cdeg}#1}
\newcommand{\w}[0]{\mathrm w}
\newcommand{\ol}[1]{\mathbf{#1}}

\title[Combinatorial Polarization]{Combinatorial Polarization, Code Loops and
Codes of High Level}

\author{Petr Vojt\v{e}chovsk\'y}

\address{Department of Mathematics, University of Denver, 2360 S Gaylord St,
Denver, CO, 80208, U.S.A.}

\email{petr@math.du.edu}

\begin{document}



\thanks{Work partially supported by Grant Agency of Charles University,
grant number $269/2001/$B-MAT/MFF}

\begin{abstract}
We first find the combinatorial degree of any map $f:V\to F$ where $F$ is a
finite field and $V$ is a finite-dimensional vector space over $F$. We then
simplify and generalize a certain construction due to Chein and Goodaire that
was used in characterizing code loops as finite Moufang loops that posses at
most two squares. The construction yields binary codes of high divisibility
level with prescribed Hamming weights of intersections of codewords.
\end{abstract}

\maketitle

\section{Introduction}

\noindent Let $V$ be a finite-dimensional vector space over $F$. Given a map
$f:V\to F$, the \emph{$n$th derived form} $\df{f}{n}:V^n\to F$ is defined by
\begin{equation}\label{Eq:DF}
    \df{f}{n}(v_1,\dots,v_n) = \sum
    (-1)^{n-m}f(v_{i_1}+\dots+v_{i_m}),
\end{equation}
where the summation runs over all subsets of $\{v_1,\dots,v_n\}$, where we
consider $v_i$ different from $v_j$ if $i\ne j$. We have borrowed the name
\emph{derived form} from \cite[p.\ 41]{Aschbacher}. The process of obtaining
the forms $\df{f}{n}$ from $f$ is referred to as \emph{combinatorial
polarization}. It can also be described recursively by the formula
\begin{multline}\label{Eq:RDF}
    \df{f}{n+1}(v_1,\dots,v_{n+1})=
    \df{f}{n}(v_1+v_2,v_3,\dots,v_{n+1})\\
    -\df{f}{n}(v_1,v_3,\dots,v_{n+1})
    -\df{f}{n}(v_2,v_3,\dots,v_{n+1}),
\end{multline}
for $v_1$, $\dots$, $v_{n+1}\in V$.

The \emph{combinatorial degree} $\cdeg{f}$ of $f$ is the smallest nonnegative
integer $n$ such that $\df{f}{m}=0$ for every $m>n$, if it exists, and it is
equal to $\infty$ otherwise. Note that the combinatorial degree of the zero map
is $0$.

Combinatorial polarization was first studied by Ward in \cite{Ward}. He shows
in \cite[Prop.\ 2.8]{Ward} that the combinatorial degree of a polynomial map is
equal to its degree if $F$ is a prime field or a field of characteristic $0$,
and he remarks in parentheses on page $195$ that ``It is not difficult to show
that, in general, the combinatorial degree of a nonzero polynomial over
$GF(q)$, $q$ a power of the prime $p$, is the largest value of the sum of the
$p$-weights of the exponents for the monomials appearing in the polynomials.''
Since he does not prove the statement and since the author believes the proof
is not absolutely trivial, we offer it in Section \ref{Sc:CombPol}. We then
move on to code loops and codes of high level.

Recall that \emph{loop} is a groupoid $(L,\cdot)$ with neutral element $1$ such
that the equation $x\cdot y=z$ has a unique solution whenever two of the three
elements $x$, $y$, $z\in L$ are specified. The variety of loops defined by the
identity $x\cdot(y\cdot(x\cdot z)) = ((x\cdot y)\cdot x)\cdot z$ is known as
\emph{Moufang loops}. The \emph{commutator} $[x,y]$ of $x$, $y\in L$ is the
unique element $u\in L$ such that $xy=(yx)u$. Similarly, the \emph{associator}
$[x,y,z]$ of $x$, $y$, $z\in L$ is the unique element $u\in L$ such that $(xy)z
= (x(yz))u$.

We shall need basic terminology from coding theory. A (linear) \emph{code} $C$
is a subspace of $V$. For every \emph{codeword} $c=(c_1,\dots,c_m)\in C$ we
define its (Hamming) \emph{weight} $\w(c)$ as the number of nonzero coordinates
$c_i$ in $c$. If $2^r$ divides $\w(c)$ for every $c\in C$ and if $r$ is as big
as possible then $C$ is said to be of \emph{level} $r$.

Codes of level $2$ are usually called \emph{doubly even}, and doubly even codes
are behind Griess' definition of a code loop: Let $C$ be a doubly even code
over $F=\{0,1\}$, and let $\eta:C\times C\to F$ be such that
\begin{eqnarray*}
    &&\eta(x,x)=\w(x)/4,\label{Eq:Griess1}\\
    &&\eta(x,y)+\eta(y,x)=\w(x\cap y)/2,\label{Eq:Griess2}\\
    &&\eta(x,y)+\eta(x+y,z)+\eta(y,z)+\eta(x,y+z) =
        \w(x\cap y\cap z)\label{Eq:Griess3}
\end{eqnarray*}
for every $x$, $y$, $z\in C$. (Here, $x\cap y$ is the vector whose $i$th
coordinate is equal to $1$ if and only if the $i$th coordinates of both $x$ and
$y$ are equal to $1$.) Then $L=C\times F$ with multiplication
\begin{equation}\label{Eq:CodeLoop}
    (x,a)(y,b)=(x+y,a+b+\eta(x,y))
\end{equation}
is a \emph{code loop} for $C$. Griess shows in \cite{Griess} that there is a
unique code loop for $C$, up to isomorphism, and that it is Moufang.

Chein and Goodaire found a nice characterization of code loops. Namely, they
show (cf. \cite[Thm.\ 5]{CheinGoodaire}) that code loops are exactly finite
Moufang loops with at most two squares. Their proof is based on three
observations:

First, if $L$ is a Moufang loop with $|L^2|\le 2$ then every commutator and
associator belongs to $L^2$ and
\begin{eqnarray}
    (xy)^2&=&x^2y^2\left[x,y\right],\notag \\
    \left[xy,z\right] &=& \left[x,z\right] \left[y,z\right]
        \left[x,y,z\right],\label{Eq:First}\\
    \left[ vx,y,z\right] &=& \left[v,y,z\right] \left[x,y,z\right]\notag
\end{eqnarray}
holds for every $v$, $x$, $y$, $z\in L$ (see \cite[Thm.\ 1,\
2]{CheinGoodaire}). In other words, if we set $Z=L^2$ then $L/Z$ is an
elementary abelian $2$-group, and the well-defined map $P:L/Z\to Z$, $x\mapsto
x^2$ satisfies $\df{P}{2}(x,y)=[x,y]$, $\df{P}{3}(x,y,z)=[x,y,z]$,
$\cdeg{P}=3$, as can be seen immediately from $(\ref{Eq:RDF})$ and
$(\ref{Eq:First})$. Note that under these circumstances $L$ is an elementary
abelian $2$-group if and only if $|L^2|=1$.

Second, if $L=C\times F$ is a code loop for $C$ and $x=(\tilde{x},a)$,
$y=(\tilde{y},b)$, $z=(\tilde{z},c)$ belong to $L$ then
\begin{eqnarray}
    x^2&=&(0,\w(\tilde{x})/4),\notag\\
    \left[x,y\right]&=&(0,\w(\tilde{x}\cap\tilde{y})/2),\label{Eq:Second}\\
    \left[x,y,z\right]&=&(0,\w(\tilde{x}\cap\tilde{y}\cap\tilde{z})),\notag
\end{eqnarray}
by \cite[Lm.\ 6]{CheinGoodaire}. (Note that $(\ref{Eq:First})$ then holds
because $\w(u+v)=\w(u)+\w(v)-2\w(u\cap v)$ for any two binary vectors $u$,
$v$.)

Third, given an integer $n\ge 1$ and parameters $\alpha_i$, $\beta_{ij}$,
$\gamma_{ijk}\in \{0,1\}$, for $1\le i$, $j$, $k\le n$, there is a doubly even
code $C$ with basis $c_1$, $\dots$, $c_n$ such that
\begin{eqnarray}
    \alpha_i&=&\w(c_i),\notag\\
    \beta_{ij}&=&\w(c_i\cap c_j),\label{Eq:CG}\\
    \gamma_{ijk}&=&\w(c_i\cap c_j\cap c_k)\notag
\end{eqnarray}
for $1\le i$, $j$, $k\le n$ (cf.\ the proof of \cite[Thm.\ 5]{CheinGoodaire}).
It is this construction that turns out to be the most difficult part of the
proof that code loops are exactly finite Moufang loops with at most two
squares. We simplify and generalize the construction in Section \ref{Sc:Code}.
The construction presented here is easier than that of \cite{Loops99}, too,
because it avoids induction.

To conclude our discussion concerning code loops, note that a map $f:V\to
\{0,1\}$ with combinatorial degree $3$ is uniquely specified if we know the
values of $f(e_i)$, $f(e_i+e_j)$ and $f(e_i+e_j+e_k)$ for some basis $e_1$,
$\dots$, $e_n$ of $V$. Hence, by $(\ref{Eq:First})$, $(\ref{Eq:Second})$ and
$(\ref{Eq:CG})$, code loops can be identified with maps $P:V\to\{0,1\}$ of
combinatorial degree $3$.

\section{Combinatorial Degree over Finite Fields}\label{Sc:CombPol}

\noindent In this section, let $V$ be a vector space of dimension $n$ over the
field $F=GF(q)$ of characteristic $p$, and let $f:V\to F$ be an arbitrary map.

By an easy generalization of the Fundamental Theorem of Algebra, we can
identify $f$ with some polynomial in $F[x_1,\dots,x_n]$. Moreover, if we assume
that the polynomial is \emph{reduced}, then it is unique. Here is what we mean
by reduced polynomial: if $cx_1^{a_1}\cdots x_n^{a_n}$ is a monomial of $f$
then $0\le a_i<q$, for $1\le i\le n$, and if $cx_1^{a_1}\cdots x_n^{a_n}$,
$dx_1^{b_1}\cdots x_n^{b_n}$ are two such monomials of $f$, we must have
$a_i\ne b_i$ for some $1\le i\le n$.

We will assume from now on that $f$ is a reduced polynomial. Let us write
$\ol{x}$ for $(x_1$, $\dots$, $x_n)$, and $\ol{x}^{\ol{a}}$ for
$x_1^{a_1}\cdots x_n^{a_n}$. Then
\begin{equation}\label{Eq:Reduced}
    f(\ol{x})=\sum_{\ol{a}\in M(f)}c_{\ol{a}}\ol{x}^\ol{a},
\end{equation}
where $M(f)$ is the set of all multiexponents of $f$, and $0\ne c_\ol{a}\in F$
for every $\ol{a}\in M(f)$.

Our goal is to calculate the combinatorial degree of $f$. When $f(0)\ne 0$ then
$\cdeg{f}=\infty$, by $(\ref{Eq:RDF})$. We will therefore assume that $f(0)=0$
from now on.

\begin{lemma}\label{Lm:FullExpansion}
Let $f:V\to F$ be as in $(\ref{Eq:Reduced})$. Then
\begin{displaymath}
    \cdeg{f} = \max\{\cdeg{\ \ol{x}^\ol{a}};\; \ol{a}\in M(f)\}.
\end{displaymath}
\end{lemma}
\begin{proof}
Let us call two polynomials $f$, $g$ disjoint if $M(f)\cap M(g)=\emptyset$.
Since every polynomial is a sum of monomials, it suffices to prove that if two
monomials are disjoint then their derived forms are disjoint, too.

To see this, consider the monomial $g(\ol{x_1})=\ol{x_1}^\ol{a}$ where
$\ol{x_1}=(x_{11}$, $\dots$, $x_{1n})$, $\ol{a}=(a_1$, $\dots$, $a_n)$. As seen
in $(\ref{Eq:DF})$, a typical summand of $\df{g}{r}$ is $h=g(\ol{x_1}+\cdots
+\ol{x_r})$, which is a polynomial in $nr$ variables. The crucial observation
is that for every $1\le j\le r$, every monomial of $h$ contains exactly $a_i$
variables $x_{ji}$, with possible repetitions. Hence the original monomial
$\ol{x_1}^\ol{a}$ can be uniquely reconstructed from every monomial of $h$.
\end{proof}

We focus on reduced monomials from now on. When $\ol{a}=(a_1$, $\dots$, $a_n)$,
$\ol{b}=(b_1$, $\dots$, $b_n)$ are two multiexponents, let us write
$\ol{a}\le\ol{b}$ if $a_i\le b_i$ for every $1\le i\le n$, and $\ol{a}<\ol{b}$
if $\ol{a}\le\ol{b}$ and there is $1\le i\le n$ with $a_i<b_i$.

\begin{lemma}\label{Lm:FirstExpansion}
Let $f(\ol{x})=\ol{x}^\ol{a}$. Then
\begin{equation}\label{Eq:Binomial}
    \df{f}{2}(\ol{x},\ol{y})=\sum_{\ol{0}<\ol{b}<\ol{a}}
        c_{\ol{a},\ol{b}}\ol{x}^\ol{b}\ol{y}^{\ol{a}-\ol{b}},
\end{equation}
where
\begin{equation}\label{Eq:ProdCoeff}
    c_{\ol{a},\ol{b}}=\prod_{i=1}^n\binom{a_i}{b_i}.
\end{equation}
More generally, if $f(\ol{x_1})=\ol{x_1}^\ol{a_1}$ then
\begin{multline}\label{Eq:FullExpansion}
    \df{f}{s}(\ol{x_1},\dots,\ol{x_s})\\
    = \sum_{\ol{0}<\ol{a_2}<\ol{a_1}} \cdots
    \sum_{\ol{0}<\ol{a_s}<\ol{a_{s-1}}}
    c_{\ol{a_1},\ol{a_2}}\cdots c_{\ol{a_{s-1}},\ol{a_s}}
    \ol{x_1}^\ol{a_s}\ol{x_2}^{\ol{a_{s-1}-\ol{a_s}}}
    \cdots \ol{x_s}^{\ol{a_1}-\ol{a_2}},
\end{multline}
where $c_{\ol{a_i},\ol{a_{i+1}}}$ is analogous to $(\ref{Eq:ProdCoeff})$.
\end{lemma}
\begin{proof}
We have
\begin{eqnarray*}
    (\ol{x}+\ol{y})^\ol{a}&=&\prod_{i=1}^n (x_i+y_i)^{a_i}
    =\prod_{i=1}^n\sum_{b_i=0}^{a_i}\binom{a_i}{b_i}x_i^{b_i}y_i^{a_i-b_i}\\
    &=&\sum_{b_1,\dots,b_n,\, 0\le b_i\le a_i}\left(
           \prod_{i=1}^n\binom{a_i}{b_i}x_i^{b_i}y_i^{a_i-b_i}\right)\\
    &=&\sum_{b_1,\dots,b_n,\, 0\le b_i\le a_i}
        \left(\prod_{i=1}^n \binom{a_i}{b_i}\right)
        x_1^{b_1}\cdots x_n^{b_n}y_1^{a_1-b_1}\cdots y_n^{a_n-b_n}.
\end{eqnarray*}
Since $\df{f}{}(\ol{x},\ol{y})=(\ol{x}+\ol{y})^\ol{a}-\ol{x}^\ol{a}
-\ol{y}^\ol{a}$, we are done with $(\ref{Eq:Binomial})$.

We have also proved the more general statement $(\ref{Eq:FullExpansion})$ for
$s=2$, and we proceed by induction on $s$. Assume that
$(\ref{Eq:FullExpansion})$ holds for $s$. Using the polarization formula
$(\ref{Eq:RDF})$ on every summand of $(\ref{Eq:FullExpansion})$, we see that
$\df{f}{s+1}(\ol{x_1},\dots,\ol{x_{s+1}})$ is equal to
\begin{displaymath}
    \sum_{\ol{0}<\ol{a_2}<\ol{a_1}}\cdots
    \sum_{\ol{0}<\ol{a_s}<\ol{a_{s-1}}}
    c_{\ol{a_1},\ol{a_2}}\cdots c_{\ol{a_{s-1}},\ol{a_s}}
    \df{g}{2}(\ol{x_1},\ol{x_2})
    \ol{x_3}^{\ol{a_{s-1}-\ol{a_s}}}
    \cdots \ol{x_{s+1}}^{\ol{a_1}-\ol{a_2}},
\end{displaymath}
where $g(\ol{x})=\ol{x}^\ol{a_s}$. The term $\df{g}{2}(\ol{x_1},\ol{x_2})$
expands as
\begin{displaymath}
     \sum_{\ol{0}<\ol{a_{s+1}}<\ol{a_s}}
        c_{\ol{a_s},\ol{a_{s+1}}}\ol{x_1}^\ol{a_{s+1}}
        \ol{x_2}^{\ol{a_s}-\ol{a_{s+1}}},
\end{displaymath}
by $(\ref{Eq:Binomial})$, and we are done.
\end{proof}

Let $\ol{a}=(a_1,\dots,a_n)$ be a multiexponent, $f(\ol{x})=\ol{x}^\ol{a}$.
Lemma \ref{Lm:FullExpansion} shows that $\df{f}{s}$ is not the zero map if and
only if there is a chain of multiexponents
$\ol{a}=\ol{a_1}>\ol{a_2}>\cdots>\ol{a_s}$ such that
$c_{\ol{a_i},\ol{a_{i+1}}}$ does not vanish in $F$ for every $1\le i<s$. We
will call such chains of multiexponents \emph{regular} here. Obviously, the
length of a regular chain is bounded by $q^n$.

\begin{lemma}
Let $\ol{a_i}=(a_{i1},\dots,a_{in})$, for $1\le i\le s$. Assume that
$\ol{a_1}>\cdots>\ol{a_s}$ is a regular chain of maximum length for $\ol{a_1}$.
Then $\ol{a_{i+1}}$, $\ol{a_i}$ differ in exactly one position, i.e.,
$a_{i+1,j}\ne a_{ij}$ for a unique $1\le j\le n$.
\end{lemma}
\begin{proof}
Suppose that there are $1\le i<s$ and $1\le j<k\le n$ such that $a_{i+1,j}\ne
a_{ij}$ and $a_{i+1,k}\ne a_{ik}$. Construct a multiexponent $\ol{b}$ so that
\begin{displaymath}
    b_m=\left\{\begin{array}{ll}
        a_{i,m},&\text{if\ }m\ne j,\\
        a_{i+1,m},&\text{if\ }m=j.
    \end{array}\right.
\end{displaymath}
Then $\ol{a_i}>\ol{b}>\ol{a_{i+1}}$. Since
$c_{\ol{a_i},\ol{a_{i+1}}}=\prod_{m=1}^n\binom{a_{i,m}}{a_{i+1,m}}\ne 0$, we
have $c_{\ol{a_i},\ol{b}}\ne 0$ and $c_{\ol{b},\ol{a_{i+1}}}\ne 0$. Thus
$\ol{a_1}$, $\dots$, $\ol{a_i}$, $\ol{b}$, $\ol{a_{i+1}}$, $\dots$, $\ol{a_s}$
is a regular chain of length $s+1$, a contradiction.
\end{proof}

\begin{corollary}
Let $\ol{a}=(a_1,\dots,a_n)$ be a multiexponent, $f(\ol{x})=\ol{x}^\ol{a}$, and
$f_i(x)=x^{a_i}$, $1\le i\le n$. Then
\begin{displaymath}
    \cdeg{f}=\sum_{i=1}^n\cdeg{f_i}.
\end{displaymath}
\end{corollary}

We therefore continue to investigate combinatorial degrees of reduced monomials
in one variable.

Let $m$ be a positive integer. Then there are uniquely determined integers
$0\le m_i<p$ such that $m=\sum_{i=0}^\infty m_ip^i$. This $p$-adic expansion of
$m$ is useful when calculating binomial coefficients modulo $p$, as seen in the
beautiful and still not so well known theorem of Lucas (cf.\ \cite{Fine}):

\begin{theorem}[Lucas Theorem]\label{Th:Lucas}
Let $p$ be a prime, and $m=\sum_{i=0}^\infty m_ip^i$, $k=\sum_{i=0}^\infty
k_ip^i$ the $p$-adic expansions of $m$ and $k$, respectively. Then
\begin{displaymath}
    \binom{m}{k}\equiv\prod_{i=0}^\infty \binom{m_i}{k_i}\pmod p,
\end{displaymath}
where we set $\binom{a}{b}=0$ whenever $a<b$.
\end{theorem}

The \emph{$p$-weight} of $m=\sum_{i=0}^\infty m_ip^i$ is defined as
$\w_p(m)=\sum_{i=0}^\infty m_i$.

\begin{corollary}
Let $p$ be a prime and $a>0$ an integer. Then the longest regular chain for $a$
has length $\w_p(a)$.
\end{corollary}
\begin{proof}
Let $\ell$ be the length of the longest regular chain for $a$. Theorem
\ref{Th:Lucas} shows that $\binom{m}{k}\not\equiv 0\pmod p$ if and only if
$m_i\ge k_i$ for every $i$ (as both $m_i$, $k_i$ are less than $p$ and
$\binom{m_i}{k_i}$ with $m_i\ge k_i$ is therefore not divisible by $p$). This
means that $\w_p(m)\ge\w_p(k)$ must be satisfied whenever
$\binom{m}{k}\not\equiv 0\pmod p$, and $\ell\le\w_p(a)$ follows.

On the other hand, if $k$ is such that $k_i=m_i$ for each $i\ne j$, and
$k_j=m_j-1\ge 0$ then $\binom{m}{k}\not\equiv 0\pmod p$, by Theorem
\ref{Th:Lucas}. Hence $\ell\ge\w_p(a)$.
\end{proof}

We have proved:

\begin{theorem}\label{Th:CD}
Let $F$ be a finite field of characteristic $p$, let $V$ be an $n$-dimensional
vector space over $F$, and let $f:V\to F$ be a map. Then $f:V\to F$ can be
written as a reduced polynomial $f(\ol{x})=\sum_{\ol{a}\in M(f)}\ol{x}^\ol{a}$
in $F[\ol{x}]$, where $\ol{x}=(x_1$, $\dots$, $x_n)$, $\ol{a}=(a_1,\dots,a_n)$,
and $M(f)$ is the set of all multiexponents of $f$. Moreover,
\begin{displaymath}
    \cdeg{f}=\left\{\begin{array}{ll}
        0,&\text{if\ }f\text{\ is the zero map},\\
        \infty,&\text{if\ }f(0)\ne 0,\\
        \deg_p{f},&\text{otherwise},
    \end{array}\right.
\end{displaymath}
where the $p$-degree $\deg_p{f}$ of $f$ is calculated as
\begin{displaymath}
       \deg_p{f}=\max_{(a_1,\dots,a_n)\in M(f)}\sum_{i=1}^n\w_p(a_i),
\end{displaymath}
where the $p$-weight $\w_p(a_i)$ of $a_i=\sum_{j=0}^\infty a_{ij}p^j$, $0\le
a_{ij}<p$, is the integer
\begin{displaymath}
    \w_p(a_i)=\sum_{j=0}^\infty a_{ij}.
\end{displaymath}
When $F$ is a prime field, we have $\deg_p f=\deg f$.
\end{theorem}

\begin{example}
Let $F=GF(9)$, $V=F^3$, $f(x_1,x_2,x_3) = f_1(x_1,x_2,x_3)+f_2(x_1,x_2,x_3)$,
where $f_1(x_1,x_2,x_3)=x_1^3x_2^7$ and $f_2(x_1,x_2,x_3)=x_1x_2x_3^5$. Since
$3=1\cdot 3^1$, $7=1\cdot 3^0+2\cdot 3^1$, $1=1\cdot 3^0$ and $5=2\cdot
3^0+1\cdot 3^1$, we have $\w_3(3)=1$, $\w_3(7)=1+2=3$, $\w_3(1)=1$ and
$\w_3(5)=2+1=3$. Thus $\cdeg{f_1}=1+3=4$, $\cdeg{f_2}=1+1+3=5$, and
$\cdeg{f}=\max\{4,5\}=5$.
\end{example}

\section{Codes of High Level with Prescribed Weights of Intersections of
Codewords}\label{Sc:Code}

\noindent We will assume from now on that $F=\{0,1\}$ is the two-element field
and that $V$ is a vector space over $F$ of dimension $n$.

In the Introduction, we have discussed a construction due to Chein and Goodaire
that is used in characterizing code loops as finite Moufang loops with at most
two squares. The fact that the values $\alpha_i$, $\beta_{ij}$, $\gamma_{ijk}$
in $(\ref{Eq:CG})$ can be prescribed can be restated as follows:

\begin{proposition}[Chein and Goodaire] Let $P:V\to F$ be a map with
$\cdeg{P}=3$. Then there is a doubly even code $C$ isomorphic to $V$ such that
$P(c)\equiv \w(c)/4\pmod 2$, for every $c\in C$.
\end{proposition}

The original proof of this Proposition is somewhat involved, and presents the
biggest obstacle when characterizing code loops. We offer a simpler proof,
while at the same time generalizing the result to codes of arbitrary level:

\begin{theorem}\label{Th:Code}
Let $V$ be an $n$-dimensional vector space over $F=\{0,1\}$, and let $P:V\to F$
be such that $\cdeg{P}=r+1$. Then there is a binary code $C$ isomorphic to $V$
and of level $r$ such that $P(c)\equiv \w(c)/2^r\pmod 2$ for every $c\in C$.
\end{theorem}
\begin{proof}
The map $P$ can be identified with some polynomial in $F[x_1,\dots,x_n]$.
Calculating in $F[x_1,\dots, x_n]$, we have
\begin{displaymath}
    \prod_{i\in I}x_i =\sum_{K\in\mathcal K}\Bigl(1+\prod_{j\in K}(1+x_j)\Bigr),
\end{displaymath}
where $I$ is a subset of $\{1$, $\dots$, $n\}$, and $\mathcal K$ is some set of
subsets of $I$. Therefore
\begin{displaymath}
    P(x_1,\dots,x_n) = \sum_{J\in\mathcal J}\Bigl(1+\prod_{j\in
    J}(1+x_j)\Bigr)
\end{displaymath}
for some set $\mathcal J$ of subsets of $\{1$, $\dots$, $n\}$.

Let $H$ be the parity-check matrix of the Hamming code of dimension $r+1$ (and
length $2^{r+1}-1$). Hence the rows of $H$ are exactly the nonzero vectors of
$F^{r+1}$, in some order. Let $D$ be the code whose generating matrix is the
transpose of $H$. Then $\w(d)=2^r$ for every nonzero $d\in D$. Every codeword
$d$ can be written as a linear combination of columns of $H$, and hence
identified with some $(d_1$, $\dots$, $d_{r+1})\in F^{r+1}$. Note that
\begin{equation}\label{Eq:W2}
    \w(d)/2^r \equiv 1+\prod_{i=1}^{r+1} (1+d_i),
\end{equation}
since the product $\prod_{i=1}^{r+1}(1+d_i)$ vanishes for every nonzero $d\in
D$.

For every subset $J=\{j_1$, $\dots$, $j_t\}$ in $\mathcal J$ define the map
$\pi_J: V\to D$ by
\begin{displaymath}
    \pi_J(x_1, \dots, x_n) = (x_{j_1}, \dots, x_{j_t}, 0, \dots, 0)\in F^{r+1}.
\end{displaymath}
This map is well-defined because $r+1=\cdeg{P}=\deg_2{P} = \deg{P} =
\max_{J\in\mathcal J}|J|$, by Theorem \ref{Th:CD}. For $x\in V$, let
$\pi(x)=\bigoplus_{J\in\mathcal J}\pi_J(x)$, and let $C$ be the image of $V$
under $\pi$. Then $C$ is isomorphic to $V$ and, for $x=(x_1$, $\dots$, $x_m)\in
V$,
\begin{displaymath}
    P(x)=
    \sum_{J\in\mathcal J}\Bigl(1+\prod_{j\in J}(1+x_j)\Bigr)
    \stackrel{(\ref{Eq:W2})}{\equiv}
    \sum_{J\in \mathcal J} \w(\pi_J(x))/2^r
     = \w(\pi(x))/2^r.
\end{displaymath}
This finishes the proof.
\end{proof}

\begin{example}
We will work out an example illustrating the proof of Theorem \ref{Th:Code}.
Let $P:V=F^3\to F$ be the map $P(x_1,x_2,x_3)=x_2+x_1x_3+x_1x_2x_3$. Then
\begin{displaymath}
    P(x_1,x_2,x_3)=(1+x_1'x_2')+(1+x_2'x_3')
        +(1+x_1'x_2'x_3'),
\end{displaymath}
where $x_i'=1+x_i$. Thus $\mathcal J=\{\{1,2\}$, $\{2,3\}$, $\{1,2,3\}\}$. We
have $\cdeg{P}=\deg{P}=3=r+1$, and $n=\dim{V}=3$.

The construction depends on a choice of the (dual) Hamming code. Let us pick
the code whose generating matrix is
\begin{displaymath}
    H^T=\left(\begin{array}{c}1000111\\0101011\\0011101\end{array}\right).
\end{displaymath}
The explicit construction also depends on an ordering of the elements of
$\mathcal J$. Let us agree that they are listed as above. Then
$\pi(x)=\pi_{\{1,2\}}(x)\oplus \pi_{\{2,3\}}(x) \oplus \pi_{\{1,2,3\}}(x)$.

Let $e_1=(1,0,0)$, $e_2=(0,1,0)$, $e_3=(0,0,1)$ be the canonical basis for $V$.
Then, with respect to the basis consisting of the three rows of $H^T$, the
vectors $e_1$, $e_2$, $e_3$ are mapped onto
\begin{eqnarray*}
    \pi(e_1)&=&(1,0,0)\oplus(0,0,0)\oplus(1,0,0)=c_1,\\
    \pi(e_2)&=&(0,1,0)\oplus(1,0,0)\oplus(0,1,0)=c_2,\\
    \pi(e_3)&=&(0,0,0)\oplus(0,1,0)\oplus(0,0,1)=c_3.
\end{eqnarray*}
For instance, the middle summand of $c_3$ is $(0,1,0)$ because
$\pi_{\{2,3\}}(0,0,1)=(0,1,0)$.

Viewed as vectors over $F$, the vectors $c_1$, $c_2$, $c_3$ (forming a basis of
$C$) become
\begin{eqnarray*}
    c_1&=&(1000111,0000000,1000111),\\
    c_2&=&(0101011,1000111,0101011),\\
    c_3&=&(0000000,0101011,0011101).
\end{eqnarray*}
Let us test Theorem \ref{Th:Code} on two vectors. First, let $x=e_3$. Then
$P(x)=0+0+0=0$, $c=\pi(x)=c_3$, and $\w(c)/2^2=8/4=2$. Similarly, with
$x=e_1+e_2+e_3$, we obtain $P(x)=1+1+1=1$,
\begin{displaymath}
    c=\pi(x)=c_1+c_2+c_3=(1101100,1101100,1110001),
\end{displaymath}
and $\w(c)/2^2=12/4=3$. In both cases, $\w(c)/2^r\equiv P(x)\pmod{2}$.
\end{example}

The construction of Theorem \ref{Th:Code} allows us to calculate the dimension
of the resulting code $C$ over $F$:

\begin{corollary}
Let $P$, $C$ be as in Theorem $\ref{Th:Code}$,
\begin{displaymath}
    P(x_1,\,\dots,\,x_n) = \sum_{J\in\mathcal J}\Bigl(1+\prod_{j\in
    J}(1+x_j)\Bigr).
\end{displaymath}
 Then the dimension of $C$ over $F$ is equal to $|\mathcal J|\cdot
 (2^{\,\deg{P}}-1)$.
\end{corollary}

\bibliographystyle{plain}

\end{document}